\documentclass[12pt]{amsart}
\usepackage[dvips]{graphicx}
\usepackage{amssymb,amsmath} 

      \newenvironment{changemargin}[2]{\begin{list}{}{
         \setlength{\topsep}{0pt}\setlength{\leftmargin}{0pt}
         \setlength{\rightmargin}{0pt}
         \setlength{\listparindent}{\parindent}
         \setlength{\itemindent}{\parindent}
         \setlength{\parsep}{0pt plus 1pt}
         \addtolength{\leftmargin}{#1}\addtolength{\rightmargin}{#2}
         }\item }{\end{list}}

\newcommand{\R}{\mathbb{R}}

\newcommand{\myfigure}[1]{\goodbreak\begin{figure}[!htp]#1\end{figure}}

\newtheorem{thm}{Theorem}[section]
\newtheorem{prop}[thm]{Proposition}
\newtheorem{cor}[thm]{Corollary}
\newtheorem{lem}[thm]{Lemma}
\newtheorem{conj}[thm]{Conjecture}

\newtheorem{notat}[thm]{Notation}
\newtheorem{alg}[thm]{Algorithm}

\theoremstyle{definition}
\newtheorem{defn}[thm]{Definition}
\theoremstyle{remark}
\newtheorem{rem}[thm]{Remark}
\newtheorem{exa}[thm]{Example}

\newcommand{\inv}{^{-1}}

\newcommand{\be}{\begin{enumerate}}
\newcommand{\ee}{\end{enumerate}}
\newcommand{\bi}{\begin{itemize}}
\newcommand{\ei}{\end{itemize}}
\renewcommand{\i}{\item}
\newcommand{\ble}{\begin{lem}}
\newcommand{\ele}{\end{lem}}
\newcommand{\bth}{\begin{thm}}
\newcommand{\bpr}{\begin{prop}}
\newcommand{\epr}{\end{prop}}
\newcommand{\bco}{\begin{cor}}
\newcommand{\eco}{\end{cor}}
\newcommand{\bcon}{\begin{conj}}
\newcommand{\econ}{\end{conj}}
\newcommand{\bde}{\begin{defn}}
\newcommand{\ede}{\end{defn}}
\newcommand{\bex}{\begin{exa}}
\newcommand{\eex}{\end{exa}}
\newcommand{\brem}{\begin{rem}}
\newcommand{\erem}{\end{rem}}
\newcommand{\bnot}{\begin{notat}}
\newcommand{\enot}{\end{notat}}
\newcommand{\balg}{\begin{alg}}
\newcommand{\ealg}{\end{alg}}

\newcommand{\N}{\mathbb{N}}

\long\def\forget#1\forgotten{} %

\newcommand{\RedGarlen}{{\ell_\mathrm{RG}}}

\begin{document}
\title[Solving equations in the braid group]{Probabilistic solutions of equations in the braid group}

\author[Garber, Kaplan, Teicher, Tsaban, Vishne]{David Garber,
  Shmuel Kaplan, Mina Teicher, Boaz Tsaban, and Uzi Vishne}

\address{David Garber,
Einstein institute of Mathematics, The Hebrew University, Givat-Ram 91904, Jerusalem, Israel; and
Department of Sciences, Holon Academic Institute of Technology, 52 Golomb Street, Holon 58102, Israel}
\email{garber@math.huji.ac.il, garber@hait.ac.il}

\address{Shmuel Kaplan, Mina Teicher, and Uzi Vishne,
Department of Mathematics and Statistics, Bar-Ilan University, Ramat-Gan 52900, Israel}
\email{[kaplansh, teicher, vishne]@math.biu.ac.il}

\address{Boaz Tsaban, Department of Applied Mathematics and Computer Science,
Weizmann Institute of Science, Rehovot 76100, Israel}
\email{boaz.tsaban@weizmann.ac.il}
\urladdr{http://www.cs.biu.ac.il/\~{}tsaban}

\thanks{\tiny This paper is a part of the Ph.D.\ thesis of the second named author
 at Bar-Ilan University.}
\thanks{\tiny
This research was partially supported by the Israel Science Foundation
through an equipment grant to the school of Computer Science in Tel-Aviv University.
The authors were partially supported by: Golda Meir
Fellowship (first named author),
EU-network HPRN-CT-2009-00099(EAGER),
Emmy Noether Research Institute for Mathematics, the Minerva
Foundation, and the Israel Science Foundation grant \#8008/02-3
(second and third named authors).}

\begin{abstract}
Given a system of equations in a ``random'' finitely generated
subgroup of the
braid group, we show how to find a small ordered list of elements
in the subgroup, which contains a solution to the equations with a
significant probability. Moreover, with a significant probability,
the solution will be the first in the list.
This gives a probabilistic solution to: The conjugacy
problem, the group membership problem, the shortest presentation
of an element, and other combinatorial group-theoretic problems in
random subgroups of the braid group.

We use a memory-based extension of the standard length-based approach,
which in principle can be applied to any group admitting an efficient,
reasonably behaving length function.
\end{abstract}

\maketitle

\section{The general method}

\subsection{Systems of equations in a group}
Fix a group $G$. A \emph{pure equation} in $G$ with variables $X_i$, $i\in\N$,
is an expression of the form
\begin{equation}\label{eq}
X_{k_1}^{\sigma_1}X_{k_2}^{\sigma_2}\dots X_{k_n}^{\sigma_n} = b,
\end{equation}
where $k_1,\dots,k_n\in\N$, $\sigma_1,\dots,\sigma_n\in\{1,-1\}$, and $b$ is given.
A \emph{parametric equation} is one obtained from a pure equation by
substituting some of the variables with given (known) parameters. By \emph{equation}
we mean either a pure or a parametric one.
Since any probabilistic method to solve a system of equations implies a probabilistic
mean to check that a given system has a solution, we will confine attention
to systems of equations which possess a solution.

Given a system of equations of the form (\ref{eq}), it is often possible
to use algebraic manipulations (taking inverses and multiplications of equations)
in order to derive from it a (possibly smaller) system of equations
all of which share the same leading variable, that is, such that all equations
have the form
\begin{equation}\label{eq2}
XW_i = b_i,
\end{equation}
where $X$ is one of the variables appearing in the original system.
The task is to find the leading variable $X$ in the system (\ref{eq2}).
Having achieved this, the process can be iterated to recover all variables
appearing in the original system (\ref{eq}).
In the sequel we confine our attention to systems
consisting of one or more equations of the form (\ref{eq2}).

\subsection{Solving equations in a finitely generated group}\label{algo}
The following general scheme is an extension of one suggested by Hughes and Tannenbaum
\cite{Hughes} and examined in \cite{Paper1}. Our new scheme turns out dramatically
more successful (compare the results of Section \ref{experiment} to those
in \cite{Paper1}).

It is convenient to think of each of the variables as an unknown element
of the group $G$.
Assume that the group $G$ is generated by the elements $a_1,\dots,a_m$,
and that there exists a ``reasonable'' length function $\ell:G\to\R^+$,
that is, such that the expected length tends to increase with
the number of multiplied generators.

Assume that equations of the form (\ref{eq2}),
$i=1,\dots,k$, are given.
We propose the following algorithm: Since $X\in G$,
it has a (shortest) form
$$X=a_{j_1}^{\sigma_1}a_{j_2}^{\sigma_2}\dots a_{j_n}^{\sigma_n}.$$
The algorithm generates an ordered list of $M$ sequences
of length $n$, such that with a significant probability, the sequence
$$((j_1,\sigma_1),(j_2,\sigma_2), \dots, (j_n,\sigma_n))$$
(which codes $X$) appears in the list, and tends to be its
\emph{first} member. The algorithm works with memory close to
$M\cdot n$, thus $M$ is usually chosen according to the memory
limitations of the computer (see also Remark \ref{complexity}).
\bi \i[\emph{Step 1}:] For each $j=1,\dots,m$ and
$\sigma\in\{1,-1\}$, compute $a_j^{-\sigma}b_i =
a_j^{-\sigma}XW_i$ for each $i=1,\dots,k$, and give $(j,\sigma)$
the score $\sum_{i=1}^k\allowbreak\ell(a_j^{-\sigma}b_i)$. Keep in
memory the $M$ elements $(j,\sigma)$ with the least scores.
\i[\emph{Step $s>1$}:] For each sequence $((j_1,\sigma_1), \dots
,(j_{s-1},\sigma_{s-1}))$ out of the $M$ sequences stored in the
memory, each $j_s=1,\dots,m$ and each $\sigma_s\in\{1,-1\}$,
compute the sum of the lengths of the elements
$$a_{j_{s}}^{-\sigma_{s}}(a_{j_{s-1}}^{-\sigma_{s-1}}\cdots a_{j_1}^{-\sigma_1}b_i) =
a_{j_{s}}^{-\sigma_{s}}a_{j_{s-1}}^{-\sigma_{s-1}}\cdots a_{j_1}^{-\sigma_1}XW_i,$$
over $i=1,\dots,k$, and assign the resulting score to the sequence
$((j_1,\sigma_1), \dots ,(j_{s},\sigma_{s}))$.
Keep in memory only the $M$ sequences with the least scores.
\ei
We still must describe the \emph{halting} condition for the algorithm.
If it is known that $X$ can be written as a product of at most $n$ generators,
then the algorithm terminates after step $n$.
Otherwise, the halting decision is more complicated. In the most general case we
can decide to stop the process when the sum of the $M$ scores
increases rather than decreases.
However, in many specific
cases the halting decision can be made much more effective --
see the examples below.

We describe several applications of the algorithm.

\begin{exa}[Parametric equations]\label{param}
If some of the words $W_i$ in the equations (\ref{eq2}) begin
with a known parameter $P_i$, then the heuristic decision when to stop
can be made much more effective: If at some step $X$ was completely
peeled of the equation, then we know the words $W_i$. To test this,
for each of the $M$ suggestions for $X$,
we calculate the words $W_i$ and check whether
the sum of the lengths $\ell(P_i\inv W_i)$ is significantly smaller than
that of the lengths $\ell(W_i)$.
In fact, this allows us to determine, with significant probability,
which of the $M$ candidates for $X$ is the correct one.
\end{exa}

\begin{exa}[The Conjugacy Problem and its variants]\label{conj}
The approach in Example \ref{param} can also be applied in the case that the system
of equations (\ref{eq2}) consists of a \emph{single} equation.
This is the case, e.g., in the \emph{parametric conjugacy problem},
where $XPX\inv$ and $P$ are given\footnote{In fact,
it is not necessary to know $P$ -- see next paragraph.}
and we wish to find $X$.
Note that in this case the algorithm can be modified to become much more
successful if at each step $s$ we peel off the generator
$a_{j_{s}}^{\sigma_{s}}$ from \emph{both sides} of the element (more precisely,
we peel off $a_{j_{s}}^{\sigma_{s}}$ from the left and $a_{j_{s}}^{-\sigma_{s}}$ from the right).

Observe, though, that if $n$ is known in advance (as in many
applications, e.g., \cite{Anshel, Ko}), then
in principle the original algorithm works, which means that we can
solve the conjugacy problem
\emph{even if we do not know the conjugated element $P$}.
\end{exa}

\begin{exa}[Group Membership and Shortest Presentation problems]\label{groupmember}
Assume that $G$ is a finitely generated subgroup of some larger group $L$.
Given $g\in L$, we wish to decide whether $g\in G$.
In this case we simply run our algorithm on $g$ using the generators of $G$,
and after each step check whether
$g$ is coded by one of our $M$ sequences.
This also provides (probabilistically) a way to write an element $g\in G$ as a product of
the generators of $G$, and with a significant
probability it will be the shortest way to write it this way.
\end{exa}

\begin{rem}[Complexity]\label{complexity}
Note that the parameter $M$ determining the length of the final list
also affects the running time of the algorithm.
As stated, if it runs $n$ steps then it performs about
$$\sum_{s=1}^nkM(s+2m)=n(n+4m+1)kM/2$$
group multiplications and $2kmnM$ evaluations of the length function
$\ell$. (Recall that $m$ denotes the number of the generators of the group,
and $k$ denotes the number of equations.)
The running time can be improved at the cost of additional memory
(e.g., one can keep in memory the $M$ elements of the form
$a_{j_{s-1}}^{-\sigma_{s-1}}\cdots a_{j_1}^{-\sigma_1}b_i$,
which were computed at step $s-1$, to reduce the number of multiplications in step $s$).
Note further that the algorithm is completely parallelable.
\end{rem}

In the next section we give experimental evidence for this algorithm's
ability to solve, with surprisingly significant probability,
arbitrary equations in ``random'' finitely generated subgroups
of the braid group $B_N$ with nontrivial parameters.

\section{Experimental results in the braid group}\label{experiment}

In the following definition (only),
we assume that the reader has some familiarity with the braid group
$B_N$ and its algorithms. Some references for these are \cite{GAR,Ko}
and references therein.

The \emph{Garside normal form} of an element $w$ in the braid group $B_N$
is a unique presentation of $w$ in the form
$\Delta_N^{-r} \cdot p_1\cdots p_m$,
where $r \geq 0$ is minimal and $p_1,\dots,p_m$ are
permutation braids in left canonical form.
The following length function was introduced in \cite{Paper1},
where it was shown that it exhibits much better properties than the
usual length function associated with the Garside normal form.
\begin{defn}[\cite{Paper1}]
Let $w=\Delta_N^{-r} \cdot p_1 \cdots p_m$ be the Garside normal
form of $w$. The \emph{Reduced Garside length} of $w$  is defined by
$$\RedGarlen(w) = r\binom{N}{2}+\sum_{i=\min\{r,m\}+1}^{m}{|p_i|}-\sum_{i=1}^{\min\{r,m\}}{|p_i|}.$$
\end{defn}

Our major experiment was made in subgroups of $B_N$ with $N=8$, which is large enough
so that $B_N$ is not trivial, but not too large so that we could perform a very large
number of experiments.
The finitely generated subgroups in which we worked were random in the sense that
each generator was chosen as a product of $10$ randomly\footnote{In this section, \emph{random}
always means with respect to the uniform distribution on the space in question.
However, we believe that good results would be obtained for any nontrivial distribution.}
chosen Artin generators.\footnote{In this section, \emph{generator} means a generator or its inverse.}
In this experiment we checked
the effectiveness of our algorithm for the \emph{parameters list}
$(m,n,k,l,M)$, where:
\be
\i $m$ (the number of generators of the subgroup) was $2$, $4$, or $8$,
\i $n$ (the number of generators multiplied to obtain $X$) was $16$, $32$, or $64$,
\i $k$ (the number of given equations of the form (\ref{eq2})) was $1$, $2$, $4$, or $8$,
\i $l$ (the number of generators multiplied to obtain the words $W_i$ in the equations (\ref{eq2})) was $4$ or $8$; and
\i $M$ (the available memory) was $2,4,8,16,32,64,128,256$, or $512$.
\ee
(see Section \ref{algo}).
This makes a total of $3\cdot 3\cdot 4\cdot 2\cdot 9=648$ parameters lists, for each
of which we repeated the experiment about $16$ times.

\subsubsection*{$X$ tends to be first}
In about $83\%$ of these experiments, $X$ was a member in the resulting list of $M$ candidates.
A natural problem is: Assume that we increase $M$.
Then experiments show that the probability of $X$ appearing in the resulting list becomes larger,\footnote{At first
glance this seems a triviality, but observe that when $M$ is increased, the correct answer has more competitors.}
but now we have more candidates for $X$, which is undesired when we cannot check which
member in the list is $X$.
However, it turns out that even for large values of $M$, $X$ tends to be among the first
few in the list.
In $71\%$ of our experiments,
$X$ was actually the first in the list, and
when $M=512$, the probabilities for $X$ ending in position $i=1,2,3,\dots$ is decreasing with $i$, and
the first few probabilities are: $0.83$, $0.08$, $0.03$, and $0.01$.

\subsubsection*{Group membership is often solved correctly}
The experiments corresponding to the group membership problem are those with $k=1$:
In these cases we are given a single element $XW$ and find a presentation of $X$
using the given generators;
this generalizes the case that we are given $X$ and find its presentation,
when it is possible (see Section \ref{groupmember}).
Checking the experiments with $k=1$, $m=4$ or $8$, and $M=512$, we get a success ratio
of $0.98$.

\subsubsection*{Logistic regression}
In order to describe the dependence of the success ratio in the
parameters involved, we are applying the methods of logistic
regression. Let $x_1,\dots,x_5$ denote the logarithms to base $2$
of the parameters $m,n,k,l,M$, respectively. Since the probability
of success $p$ in each case is a number between $0$ and $1$, a
standard linear model (expressing $p$ as a linear combination of
the variables $x_i$) is not suitable. Instead, it is customary to
express the function $L=\log(p/(1-p))$ as such a linear
combination of the variables $x_i$ (so that $p = e^L/(1+e^L)$).
This is called the \emph{logistic model}. Note that
under this transformation the derivative of $p$ with respect to
$L$ is $p(1-p)$, so an addition of $\Delta L$ to $L$ will increase
$p$ to approximately $p + p(1-p)\Delta L$. The best approximation
in this model is
\begin{equation}\label{firstL}
L \approx 7.0814 - 1.7165 x_1 - 0.7547 x_2 + 0.1094 x_3 + 0.5437 x_5.
\end{equation}
The quality of the approximation is measured by the variance
of the error. Since we are taking the best linear
approximation, adding \emph{any} variable (even a random
independent one) reduces the variance of the error. The
significance level of a variable $x_i$ roughly measures the
probability that adding this variable to the others will have its
reducing effect, assuming it was random. The typical threshold is
$0.05$: A significance level of $0.05$ or below means that the variable
has a significant contribution to the approximation $L$, which could not
be attained by a variable independent of $L$.
In the approximation (\ref{firstL}),
all variables have significance level $< 0.0003$, except for the variable
$x_4$ (corresponding to $l$) which has significance level $0.096$,
and is therefore not taken into consideration in the approximation (\ref{firstL}).

We have verified that Approximation (\ref{firstL})
gives a fairly good estimation of the success probabilities
for the tried parameters.
\subsubsection*{Doubling the memory}
Figure \ref{MultiMs} shows the effect of doubling $M$ on the success probability,
according to Approximation (\ref{firstL}).
To create this figure, we fixed $m=8$ and $k=1$, and for each $M=2^1,2^2,\dots,2^{10}$
we have drawn the graph of the success probability $p$ with respect to $\log_2(n)$.

\myfigure{
\begin{changemargin}{-2cm}{-2cm}
\begin{center}
\includegraphics[scale=0.36]{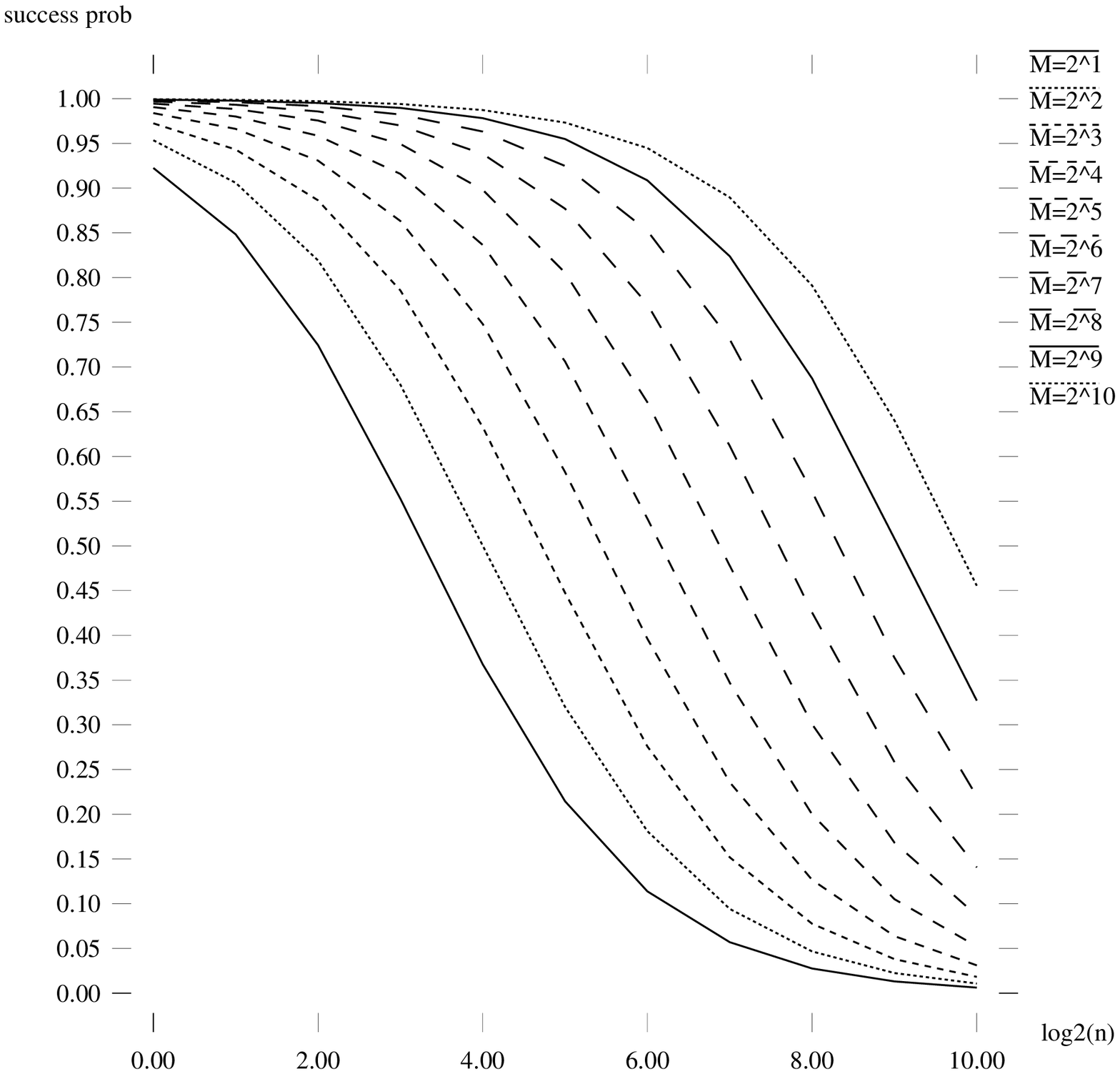}
\caption{The effect of doubling $M$ on the success probability}\label{MultiMs}
\end{center}
\end{changemargin}
}

\begin{rem}
According to Approximation (\ref{firstL}),
in order to maintain the success probability
when $m$ is doubled, $M$ should be multiplied by
$2^{1.7165/0.5437}\approx 8.92$.

Another interpretation is as follows.
Assume that we wish to decide what should the value of $M$ be to get success probability $0.5$,
that is, $L=0$. From (\ref{firstL}) it follows that
$$x_5 \approx (-7.0814 + 1.7165 x_1 + 0.7547 x_2 - 0.1094 x_3)/0.5437$$
and therefore
$$M = 2^{x_5} \approx 0.00012\cdot m^{3.16}\cdot n^{1.39}/k^{0.2}.$$
\end{rem}

It seems that the prediction capabilities of Approximation (\ref{firstL})
for larger parameters are not bad.

\begin{exa}
Using Approximation (\ref{firstL}), the predicted success probability for
parameters list $(16,128,8,8,1024)$ is $0.668$.
An experiment for these parameters succeeded in $9$ out of $11$ tries (about $0.82$).
\end{exa}

\subsection{Identifying failures}
Figure \ref{Position of correct result} describes the position of the correct prefix of $X$
and the average score of all $M$ sequences in the memory
during the steps of the algorithm (The
graphs are normalized for graphical clarity).
Two typical examples are given,
both for parameters list $(2,64,8,8,128)$.
An interesting observation is that when the correct prefix is not among the first few,
the average length decreases more slowly with the steps of the algorithms.
\myfigure{
\begin{changemargin}{-2cm}{-2cm}
\begin{center}
\includegraphics[scale=0.33]{Good2.ps}
\includegraphics[scale=0.33]{Good1.ps}\\
\caption{Position of the correct prefix in successful runs}\label{Position of correct result}
\end{center}
\end{changemargin}
}

It turns out that in most of the cases where the correct prefix of $X$ does not survive a
certain step (that is, it is not ranked among the first $M$ sequences),
the average length after several more steps almost does not decrease.
Figure \ref{failures} illustrates two typical cases, with parameters list $(2,64,8,8,16)$
(left) and $(2,64,8,8,8)$ (right).

\myfigure{
\begin{changemargin}{-2cm}{-2cm}
\begin{center}
\includegraphics[scale=0.33]{Bad1.ps}
\includegraphics[scale=0.33]{Bad2.ps}\\
\caption{Position of the correct prefix in unsuccessful runs}\label{failures}
\end{center}
\end{changemargin}
}

This allows us to identify failures within several steps after their occurrence. In such cases
one approach is to return a few steps backwards, increase $M$ for the next (problematic) few
steps, and then decrease it again.

We must stress that these are only typical cases, and several pathological cases (where the correlation
between the decrease in the lengths and the position of the correct prefix was not as expected)
were also encountered.
In these rare cases, we observed at least one of the following phenomena:
Either the generators $a_i$ could be written as a product of very
few Artin generators, due to several cancellations in the product defining them
(recall that each generator $a_i$ is a product of $10$ random Artin generators in $B_8$),
or else some (but not all) of the Artin generators
multiplied to obtain $a_i$ were cancelled when multiplied with some of the
Artin generators defining $a_j$ (or its inverse),
so that the resulting element $x$ could be written using much fewer
Artin generators than expected.
This violates the required monotonicity of the length function
and makes the algorithm fail.

\subsection{Working in $B_N$ when $N$ is larger}
For the parameters lists $(2,16,8,8,2)$ and $(8,16,8,8,128)$, we have checked the success probabilities
for $N=8$, $10$, $12$, $14$, $16$, $20$, $24$, $28$, $32$, $36$, $40$, $50$, $60$, $70$, $80$, $96$, and $100$.
The results are shown in Figure \ref{MultiBs}.
While the success probability decreases with $N$, it does not become as
negligible as one might expect. Moreover, it
can be significantly enlarged at the cost of increasing $M$.

\myfigure{
\begin{changemargin}{-2cm}{-2cm}
\begin{center}
\includegraphics[scale=0.4]{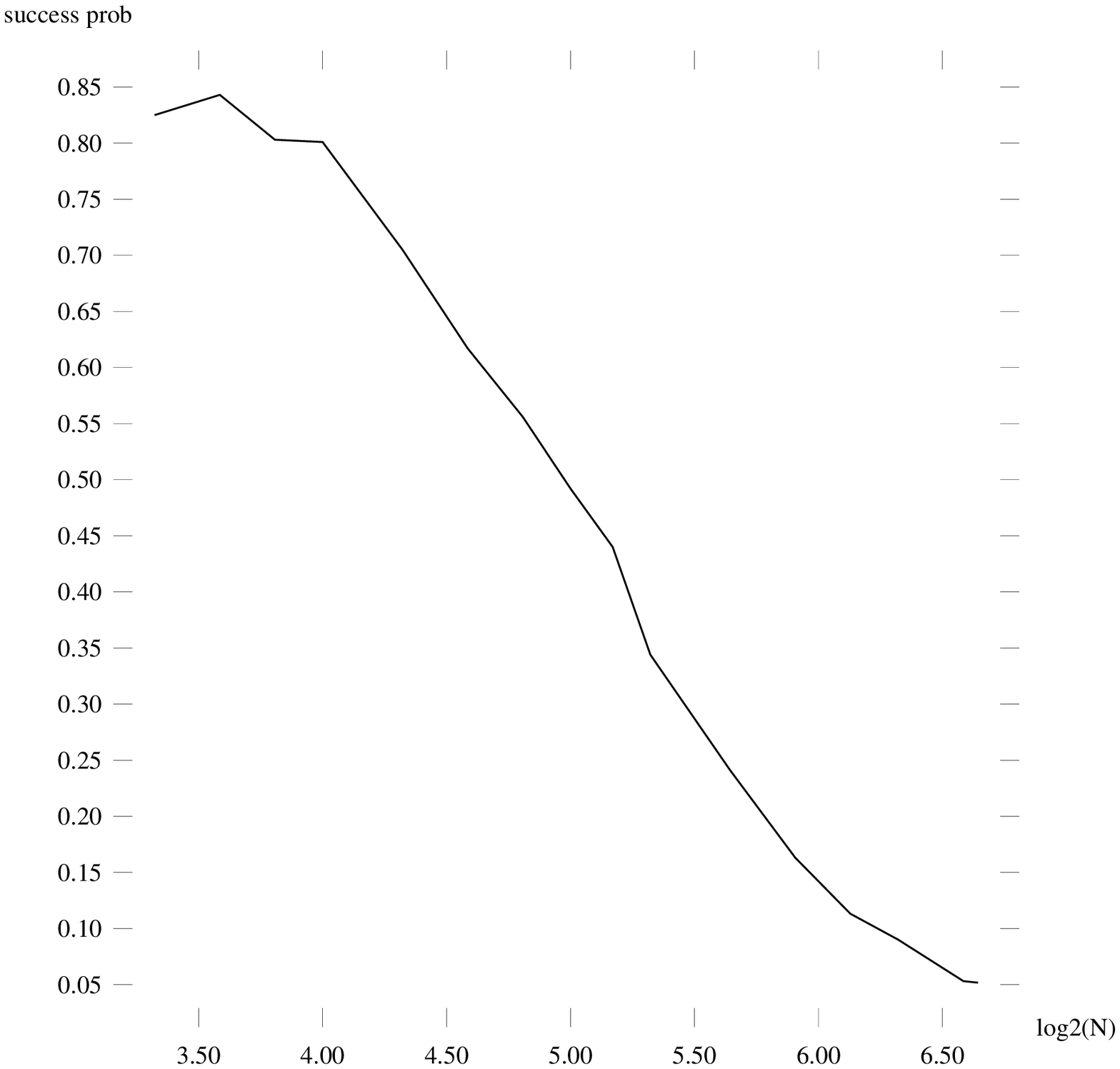}
\includegraphics[scale=0.4]{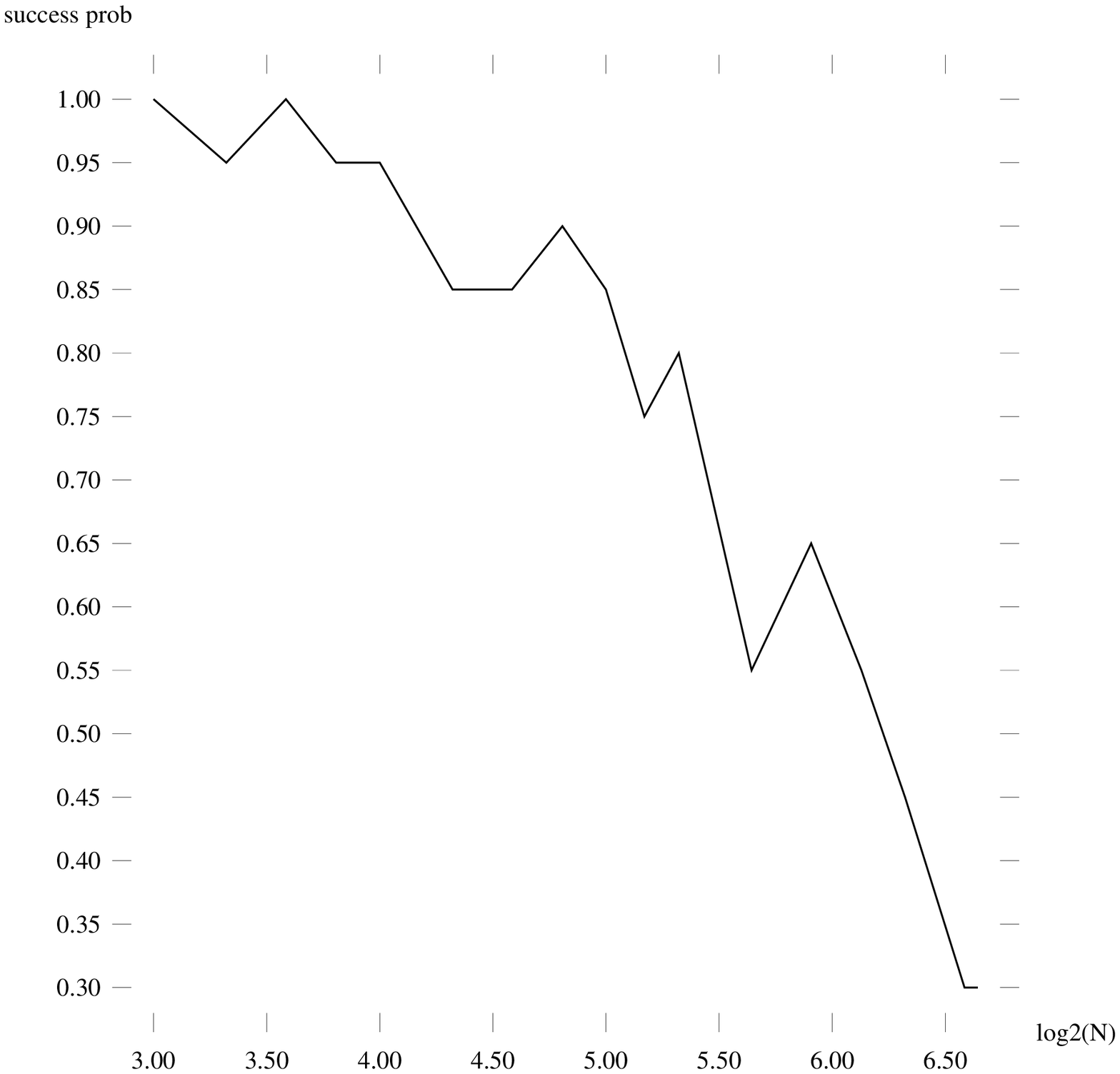}\\
\caption{Success probability for $(2,16,8,8,2)$ (left) and for $(8,16,8,8,128)$ (right)}\label{MultiBs}
\end{center}
\end{changemargin}
}

\section{Concluding remarks}
Our results suggest that whenever $G$
is a finitely generated subgroup of the braid group, which
is obtained by a sufficiently ``random'' process,
and the involved parameters are feasible for handling the group
elements in the computer, it is possible to solve equations
in the given group with significant success probabilities.
This significantly extends similar results
concerning the conjugacy problem (with known parameters) obtained
in other works (e.g., \cite{HofSte}).

This approach seems to imply the vulnerability of the
key exchange protocols suggested in \cite{Anshel, Ko},
since their security is based on the difficulty of the Conjugacy Problem
in ``random'' subgroups of the braid group (see Example \ref{conj}).
It should be stressed that our experiments were
performed with a small amount of memory (parameter $M$), which
could, in feasible settings, be increased by several orders of magnitude
and therefore significantly improve the success probability.
Since even a small non-negligible success probability in
attacking the protocol implies that it is not secure,
it seems
that in order to immune the current
protocols against the attack implied by the results here, the working
parameters have to be increased so much that the system will
become impractical.

However, in order to use
our approach against newly proposed protocols based on the braid
group (see \cite{BGC}), or against similar protocols based on
other finitely generated groups,
one must first find a good length function
for the specific problem.

\end{document}